\documentclass{lms}

\usepackage{amssymb}
\usepackage{amsmath}

\newtheorem{theorem}{Theorem}[section] 

\newtheorem{proposition}[theorem]{Proposition}

\newnumbered{assertion}{Assertion}    
\newnumbered{conjecture}{Conjecture}  
\newnumbered{definition}[theorem]{Definition}
\newnumbered{hypothesis}{Hypothesis}
\newnumbered{remark}[theorem]{Remark}
\newnumbered{note}{Note}
\newnumbered{observation}{Observation}
\newnumbered{problem}{Problem}
\newnumbered{question}{Question}
\newnumbered{algorithm}{Algorithm}
\newnumbered{example}[theorem]{Example}
\newunnumbered{notation}{Notation} 

\providecommand{\abs}[1]{\lvert#1\rvert}
\providecommand{\norm}[1]{\lVert#1\rVert}

\title[Strassen's and Choquet theorems and transport]
{Corrigendum to `Applications of Strassen's theorem and Choquet theory to optimal transport problems, to uniformly convex functions and to uniformly smooth functions'} 

\author{Krzysztof J. Ciosmak}

\dedication{}

\classno{Primary 46A55, 49N15, 26B25; Secondary 60G42, 90C46, 49N05, 47H05Declaration of interests: none. \\ This research was supported by the ERC Starting Grant 802689 CURVATURE}

\begin{document}
\maketitle

\begin{abstract}
In [K.J. Ciosmak, Applications of Strassen's theorem and Choquet theory to optimal transport problems, to uniformly convex functions and to uniformly smooth functions, Nonlinear Anal. 232 (2023), Paper No. 113267, 32 pp.], Theorem 2.3. does not suffice for its applications. We strengthen Theorem 2.1. and Theorem 2.3., so that they imply their claimed consequences. 
Moreover, we correct a minor flaw in the proof of Proposition 2.4..
\end{abstract}

\section{Notation}

Before we proceed to the corrected formulations and proofs of \cite[Theorem 2.1., Theorem 2.3.]{Ciosmak0}, let us recall a definition and notation from \cite{Ciosmak0}.

\begin{definition}
If $(\Omega,\Sigma)$ and $(\Xi, \Theta)$ are measurable spaces, then a Markov kernel $P$ from $\Omega$ to $\Xi$ is a real function on $\Theta\times \Omega$ such that for any point $\omega\in\Omega$, $P(\cdot,\omega)$ is a probability measure on $\Theta$ and for any $A\in \Theta$, $P(A,\cdot)$ is $\Sigma$-measurable. 

If $\mu$ is a probability measure on $\Sigma$, then we define $P\mu$ to be a probability measure on $\Theta$ such that
\begin{equation*}
P\mu(A)=\int_{\Omega} P(A,\omega)\,d\mu(\omega)\text{ for all }A\in\Theta.
\end{equation*}
\end{definition}

We shall denote by $\mathcal{C}(\Omega)$ the Banach space of bounded continuous functions on a topological space $\Omega$ and by $\mathcal{M}(\Omega)$ the Banach space of signed Borel measures on $\Omega$ normed by total variation. By $\mathcal{P}(\Omega)$ we shall denote the set of Borel probability measures on $\Omega$. A subset $\mathcal{K}$ of $\mathcal{C}(\Omega)$ is said to be stable under maxima provided that the maximum $f\vee g$ of any two functions   $f,g\in\mathcal{K}$ belongs to $\mathcal{K}$.

\section{New version of Strassen's theorem}

We shall now provide new versions of \cite[Theorem 2.1.]{Ciosmak0} and \cite[Theorem 2.3.]{Ciosmak0}. We discuss the introduced amendments after the respective formulations of the new versions.

\begin{theorem}\label{thm:varstrassen}
Let $X$ be a separable Banach space, let $(\Omega,\Sigma,\mu)$ be a probability space. Let $\omega\mapsto h_{\omega}$ be a map from $\Omega$ to continuous, convex functions on $X$, which is weakly measurable, that is, for every $x\in X$ the map $\omega\mapsto h_{\omega}(x)$ is $\Sigma$-measurable, and such that there exists non-negative number $c$ such that
\begin{equation}\label{eqn:ccc}
\abs{h_{\omega}(x)}\leq c(\norm{x}+1)\text{ for all }x\in X.
\end{equation}
Set 
\begin{equation*}
h(x)=\int_{\Omega}h_{\omega}(x)\,d\mu(\omega).
\end{equation*}
For a functional $x^*\in X^*$ the following conditions are equivalent:
\begin{enumerate}[i)]
\item\label{i:first} $x^*\leq h$,
\item\label{i:second} there exists a map $\omega\mapsto x_{\omega}^*$ from $\Omega$ to $X^*$ which is weakly measurable, in the sense that $\omega\mapsto x_{\omega}^*(x)$ is measurable for any $x\in X$, and such that $x_{\omega}^*\leq h_{\omega}$ for $\mu$-almost every $\omega\in\Omega$ and for all $x\in X$
\begin{equation*}
x^*(x)=\int_{\Omega}x_{\omega}^*(x)\,d\mu(\omega).
\end{equation*}
\end{enumerate}
\end{theorem}

In the theorem above, we assume that the map $\omega\mapsto h_{\omega} $ has values in continuous and convex functions, as opposed to having values in continuous, sublinear and positively homogeneous functions, which was assumed in \cite[Theorem 2.1.]{Ciosmak0}. Moreover, the assumption (\ref{eqn:ccc}) is a relaxation of the corresponding assumption of \cite[Theorem 2.1.]{Ciosmak0}.

\begin{proof}
Observe that clearly \ref{i:second} implies \ref{i:first}. For the proof of the other implication, consider the Banach space $L^1(\Omega,X)$ of equivalence classes of $\mu$-measurable, Bochner-integrable maps $f\colon \Omega\to X$ with norm
\begin{equation*}
\norm{f}_1=\int_{\Omega}\norm{f(\omega)}\,d\mu(\omega).
\end{equation*}
Define a functional $\lambda$ on the subspace of $L^1(\Omega,X)$ consisting of constant functions by the formula
\begin{equation*}
\lambda(x)=x^*(x).
\end{equation*}
By the assumption, $\lambda\leq h$, on that subspace.
By the Hahn--Banach theorem for convex majorants we may extend it to a functional $\Lambda$ defined on $L^1(\Omega, X)$ such that it satisfies
\begin{equation*}
\Lambda(g)\leq \int_{\Omega}h_{\omega}(g(\omega))\,d\mu(\omega).
\end{equation*}
Note that the right-hand side of the above inequality is well defined by the assumption (\ref{eqn:ccc}).
Then there exists a weakly measurable map $x_{\omega}^*$ such that 
\begin{equation*}
\Lambda(g)=\int_{\Omega}x_{\omega}^*(g(\omega))\,d\mu(\omega),
\end{equation*} 
cf. \cite[Theorem T51]{Meyer}.
To complete the proof of the theorem we need to show that for $\mu$-almost every $\omega$ there is $x_{\omega}^*\leq h_{\omega}$. Choose a dense subset $(x_n)_{n\in\mathbb{N}}$ in $X$. It is enough to show that for all $n\in\mathbb{N}$ there is $x_{\omega}^*(x_n)\leq h_{\omega}(x_n)$ for $\mu$-almost every $\omega$.
For $n\in\mathbb{N}$ and $A\in\Sigma$ consider map $x_n\mathbf{1}_A\in L^1(\Omega,X)$. Then it follows that
\begin{equation}\label{eqn:ae}
\int_A x_{\omega}^*(x_n)\,d\mu(\omega)\leq \int_A h_{\omega}(x_n)\,d\mu(\omega).
\end{equation}
The assertion of the theorem follows by standard arguments.
\end{proof}

\begin{theorem}\label{thm:extr}
Let $\Omega$ be a locally compact Polish space. Let $\mathcal{K}$ be a convex set in $\mathcal{C}(\Omega)$ that is stable under maxima, contains constants, and for any constant $c$ it holds $\mathcal{K}+c\subset\mathcal{K}$.
Let $f\in\mathcal{C}(\Omega)$. Suppose that $\mu,\nu$ are Borel probability measures such that 
\begin{equation}\label{eqn:maj}
\int_{\Omega}g\,d(\mu-\nu)\leq\int_{\Omega} f\,d(\mu-\nu)
\end{equation}
for all $g\in\mathcal{K}$.
Then there exists a Markov kernel $P$ form $\Omega$ to $\Omega$ such that  $\nu=P\mu$ and such that for every $\omega\in\Omega$
\begin{equation*}
\int_{\Omega}g\,d(\delta_{\omega}-P(\cdot,\omega))\leq\int_{\Omega} f\,d(\delta_{\omega}-P(\cdot,\omega))
\end{equation*}
for all $g\in\mathcal{K}$.
Moreover, the set of extreme points of the set of pairs of Borel probability measures $(\mu,\nu)$ that satisfy (\ref{eqn:maj}) is contained in the set of pairs of the form $(\delta_{\omega},\eta)$ for some $\omega\in \Omega$ and some Borel probability measure $\eta$ on $\Omega$.
\end{theorem}

The difference between the formulation above and the formulation of \cite[Theorem 2.3.]{Ciosmak0} is that above $f$ is not required to belong to $\mathcal{K}$. As follows by \cite[Proposition 2.4.]{Ciosmak0}, see the corrected Proposition \ref{pro:extr} below, if $f\in\mathcal{K}$, then the tangent cone $\mathcal{F}=\{\lambda(g-f)\mid g\in\mathcal{K},\lambda\geq 0\}$ is stable under maxima. In such case \cite[Theorem 2.3.]{Ciosmak0} is indeed sufficient. However, in \cite{Ciosmak0} at several places, in \cite[Corollary 2.6.]{Ciosmak0} and consequently in \cite[Lemma 4.2., Theorem 5.3., Theorem 8.1., Lemma 10.5.]{Ciosmak0}, stronger Theorem \ref{thm:varstrassen} is needed.

\begin{proof}
Set $X=\mathcal{C}(\Omega)$ to be the Banach space of all continuous bounded functions on $\Omega$. Let $x^*$ be an element of $X^*$ represented by a measure $\nu\in\mathcal{M}(\Omega)$. Set for $\omega\in\Omega$
\begin{equation*}
h_{\omega}(x)=\inf\big\{(f-y)(\omega)\mid y\in\mathcal{K}, f-y\geq x\big\}.
\end{equation*}
Then, as $\mathcal{K}$ is convex, $h_{\omega}$ is convex. It is moreover continuous. Indeed, choose $x_1,x_2\in \mathcal{C}(\Omega)$ and let $\delta=\norm{x_1-x_2}$.  Take $\epsilon>0$, $\omega\in\Omega$ and $y_1\in \mathcal{K}$ such that
\begin{equation*}
(f-y_1)(\omega)<h_{\omega}(x_1)+\epsilon\text{ and }x_1\leq f-y_1.
\end{equation*}
Then $y_2=y_1-\delta$ belongs to $\mathcal{K}$ and satisfies
\begin{equation*}
x_2\leq x_1+\delta\leq f-y_2\text{ and }h_{\omega}(x_2)\leq (f-y_2)(\omega)<h_{\omega}(x_1)+\delta+\epsilon.
\end{equation*}
It follows that $h_{\omega}(x_2)\leq h_{\omega}(x_1)+\delta$. By symmetry, it follows that $h_{\omega}$ is $1$-Lipschitz.
As $\mathcal{C}(\Omega)$ is separable, so is its subset 
\begin{equation*}
\big\{f-y\mid y\in\mathcal{K}, f-y\geq x\big\}.
\end{equation*}
By the assumption that $\mathcal{K}$ is stable under maxima, $\omega\mapsto h_{\omega}(x)$ is a pointwise limit of a sequence $(f-y_k)_{k=1}^{\infty}$ with $y_k\in\mathcal{K}$. We may moreover assume that 
\begin{equation*}
y_k(\omega)\geq -\norm{f-x} \text{ for }\omega\in \Omega.
\end{equation*}
Hence $(f-y_k)(\omega)\leq\norm{f-x}+\norm{f}$ for $\omega\in\Omega$. By the definition of $h_{\omega}(x)$ it follows that 
\begin{equation*}
    (f-y_k)(\omega)\geq - \norm{x}.
\end{equation*}
By the assumption on $\mu,\nu$
\begin{equation*}
x^*(x)=\int_{\Omega} x\,d\nu \leq \int_{\Omega} (f-y_k)(\omega)\,d\nu(\omega)\leq\int_{\Omega}(f-y_k)(\omega)\,d\mu(\omega).
\end{equation*}
Now, by the dominated convergence theorem it follows that
\begin{equation*}
x^*(x)\leq \int_{\Omega}h_{\omega}(x)\,d\mu(\omega).
\end{equation*} 
By the previous observations, $\abs{h_{\omega}(x)}\leq 2\norm{f}+\norm{x}$. 
Observe that the Banach space $\mathcal{C}_0(\Omega)$ is separable.  By \cite[Theorem 2.1.]{Ciosmak0}
we see that there is a weakly measurable function $\omega\mapsto x_{\omega}^*$ with values in $X^*$ such that 
\begin{equation}\label{eqn:rep}
x^*(x)=\int_{\Omega}x_{\omega}^*(x) \,d\mu(\omega)\text{ for all }x\in \mathcal{C}(\Omega)
\end{equation}
and $x_{\omega}'\leq h_{\omega}$ for $\mu$-almost every $\omega\in\Omega$.
Here $x_{\omega}'$ is the restriction of $x_{\omega}^*$ to $\mathcal{C}_0(\Omega)$.

Observe that $h_{\omega}(x)\leq\norm{f-x}+\norm{f}$. Thus
\begin{equation*}
x_{\omega}'(x)\leq h_{\omega}(x)\text{ for all }x\in \mathcal{C}_0(\Omega)
\end{equation*}
implies that for all positive $\lambda$, all $x\in\mathcal{C}_0(\Omega)$
\begin{equation*}
x_{\omega}'(x)\leq \Big\lVert \frac{f}{\lambda}-x\Big\rVert+\frac{\norm{f}}{\lambda}
\end{equation*}
Letting $\lambda$ to infinity we infer that the norm of $x_{\omega}'$ is at most one. 
By Riesz' representation theorem, there exists a measure $P(\cdot,\omega)$ on $\Omega$ such that for all $h\in \mathcal{C}_0(\Omega)$ there is 
\begin{equation*}
x_{\omega}'(h)=\int_{\Omega}h\, dP(\cdot,\omega).
\end{equation*}
Choose any $h\in \mathcal{C}(\Omega)$. By Ulam's lemma and Urysohn's lemma there exists a bounded, monotone, sequence of non-negative continuous and compactly supported functions $(\phi_n)_{n\in \mathbb{N}}$ that converges pointwise to constant function $1$. In view of (\ref{eqn:rep}),
\begin{equation*}
\int_{\Omega}h\phi_n\,d\nu=\int_{\Omega}\int_{\Omega} h\phi_n \,dP(\cdot,\omega) \,d\mu(\omega).
\end{equation*}
By the dominated convergence theorem there is 
\begin{equation}\label{eqn:repp}
\int_{\Omega}h\,d\nu=\int_{\Omega}\int_{\Omega} h \,dP(\cdot,\omega)\, d\mu(\omega).
\end{equation}
It follows that $P(\cdot,\omega)$ is a probability measure for $\mu$-almost every $\omega$. 
Since $x_{\omega}'\leq h_{\omega}$ for $\mu$-almost every $\omega\in\Omega$, we see that for all $h\in\mathcal{C}_0(\Omega)$ there is
\begin{equation*}
\int_{\Omega} h\,dP(\cdot,\omega)\leq h_{\omega}(h)\text{ for }\mu\text{-almost every }\omega\in\Omega.
\end{equation*}
As $h_{\omega}$ is monotone on $\mathcal{C}(\Omega)$, we see that the above inequality holds true for $h\in\mathcal{C}(\Omega)$ as well.

Observe that if $g\in f-\mathcal{K}$, then by the definition of $h_{\omega}$,
\begin{equation*}
h_{\omega}(g)=g\text{, hence }
\int_{\Omega} h\,dP(\cdot,\omega)\leq  g(\omega).
\end{equation*}
It follows that for all $h\in\mathcal{K}$ there is
\begin{equation}\label{eqn:ex}
\int_{\Omega} h\,d\delta_{\omega}-\int_{\Omega} h\,dP(\cdot,\omega)\leq \int_{\Omega} f\,d\delta_{\omega}-\int_{\Omega} f\,dP(\cdot,\omega).
\end{equation}

We see that $P$ defines a Markov kernel from $\Omega$ to $\Omega$. 
By (\ref{eqn:repp}), $\nu=P\mu$ and by (\ref{eqn:ex}) we have
\begin{equation}\label{eqn:modify}
\int_{\Omega} g \,d(\delta_{\omega}-P(\cdot,\omega))\leq \int_{\Omega} f \,d(\delta_{\omega}-P(\cdot,\omega))
\end{equation}
for $\mu$-almost every $\omega\in\Omega$ and all $g\in\mathcal{K}$. To obtain the desired Markov kernel we modify $P$ on a measurable set of $\mu$-measure zero of $\omega\in \Omega$ such that the inequality (\ref{eqn:modify}) is not valid for some $g\in\mathcal{K}$ by putting $P(\cdot,\omega)=\delta_{\omega}$.  

We have
\begin{equation*}
(\mu,\nu)=\int_{\Omega} (\delta_{\omega},P(\cdot,\omega))\,d\mu(\omega),
\end{equation*}
so, by (\ref{eqn:ex}), the claim about the extreme points follows.
\end{proof}

In the proposition below, we correct the erroneous definition of the tangent cone of \cite[Proposition 2.4.]{Ciosmak0} and, consequently, we provide a corrected proof.

\begin{proposition}\label{pro:extr}
Suppose that $\mathcal{K}$ is a convex set that is stable under maxima, contains constants and for any constant $c$ it holds $\mathcal{K}+c\subset\mathcal{K}$. 
Let $f\in\mathcal{K}$. Then the tangent cone to $\mathcal{K}$ at $f$, i.e., the cone
\begin{equation*}
\mathcal{F}=\Big\{\lambda(g-f)\mid g\in\mathcal{K},\lambda\geq 0\Big\}
\end{equation*}
is a convex cone, stable under maxima and contains constants.
\end{proposition}
\begin{proof}
Let $\lambda_1,\lambda_2\geq 0$, $g_1,g_2\in\mathcal{K}$. Set $\lambda=\lambda_1\vee\lambda_2$ and 
\begin{equation*}
g=\Bigg(\frac{\lambda_1}{\lambda} g_1+ \Big(1-\frac{\lambda_1}{\lambda}\Big)f \Bigg)\vee\Bigg(\frac{\lambda_2}{\lambda} g_2+ \Big(1-\frac{\lambda_2}{\lambda}\Big)f \Bigg).
\end{equation*}
Then $g\in\mathcal{K}$, thanks to the assumption that $f\in\mathcal{K}$, convexity of $\mathcal{K}$, and stability of $\mathcal{K}$ with respect to maxima.
Moreover
\begin{equation}
\Big(\lambda_1(g_1-f)\Big)\vee\Big(\lambda_2(g_2-f)\Big)=\lambda(g-f).
\end{equation}
This shows that $\mathcal{F}$ is stable under maxima. The other two claimed properties of $\mathcal{F}$ are trivial to verify.
\end{proof}


  \bibliographystyle{plain} 
  \bibliography{refs}
  
\affiliationone{
   Krzysztof J. Ciosmak\\
Fields Institute for Research in Mathematical Sciences, 222 College Street, Toronto, Ontario M5T 3J1, Canada\\
Department of Mathematics, University of Toronto, Bahen Centre, 40 St. George St., Room 6290, Toronto, Ontario, M5S 2E4, Canada
   \email{k.ciosmak@utoronto.ca}}

%


\end{document}